\documentstyle[12pt,amssymb,pb-diagram]{amsart}
\advance \topmargin by -\headheight
\advance \topmargin by -\headsep
\newtheorem{prop}{Proposition}[section]
\newtheorem{lemma}{Lemma}[section]
\newtheorem{cor}{Corollary}[section]
\title{Refining the Abel--Jacobi maps}
\author{M.Rovinsky}
\begin{document}
\begin{abstract}
Given a smooth projective variety $X$ over a field $k$ of characteristic 
zero, we consider the composition of the de Rham cohomology cycle class map 
over $k$ from the Chow group $CH^q(X\times_kK)$, where $K$ is the field of 
fractions of henselization $A^h$ of the local ring of a smooth closed 
point of a variety over the field $k$ with an appropriate projection:
$$CH^q(X\times_kK)\longrightarrow\bigoplus_{p=1}^qgr_F^{q-p}N^{q-p}
H^{2q-p}_{dR/k}(X)\otimes_k\Omega^p_{A^h/k,{\rm closed}},$$
where $F^{\bullet}$ and $N^{\bullet}$ are the Hodge and the
coniveau filtrations on the de Rham cohomology, respectively. 
The classical Abel--Jacobi map corresponds to the composition 
of this homomorphism with the projection to the summand $p=1$. 

This homomorphism is not injective, however, its composition with 
the embedding into the space 
$$\bigoplus_{p=1}^qgr_F^{q-p}N^{q-p}H^{2q-p}_{dR/k}(X)\otimes_k
\lim_{\longleftarrow_M}d(\Omega^{p-1}_{A_M/k}),$$ 
where $A_M=A^h/{\frak m}^M$ and ${\frak m}$ is the maximal 
ideal, is dominant for any $q$ for which the inverse Lefschetz operator 
$H^{2\dim X-q}(X)(\dim X)\stackrel{\sim}{\longrightarrow}H^q(X)(q)$
is induced by a correspondence. 
\end{abstract}
\maketitle
As it is shown by Mumford, if there is a holomorphic 2-form on 
a smooth complex projective surface $S$ then the group of 0-cycles 
of degree 0 on $S$ is not presentable as a quotient of an 
algebraic group, ``the group $CH_0(S)^0$ is infinite-dimensional''. 
It is clear from the argument that the ``arithmetics'' of the
complex numbers is pretty much involved. 

Given, say, a complex variety and having in mind the Mumford's 
obstruction to describing the Chow groups in geometrical terms 
and that the arithmetics of the huge field of complex numbers is 
deeply hidden, we replace the field of the complex numbers with 
the field of fractions of henselization of the local ring of a 
smooth closed point of a variety over a ``small'' subfield of ${\Bbb C}$. 

As it is well-known, the Picard group of a smooth projective
variety is an extension of a finitely generated group by an
abelian variety. We remark that the group of points of the
abelian variety over the field of fractions of an algebra of 
formal power series is an extension of the ``small'' group 
of points of the abelian variety over the field of constants 
by the maximal ideal in the algebra of formal power series 
tensored with the tangent space of the abelian variety, i.e., the 
$H^1(~\cdot~,{\cal O})$-group of the original variety. 

In general, suppose $X$ is a smooth projective variety over a field 
$k$ of characteristic zero and $K$ is the field of fractions of 
henselization $A^h$ of the local ring of a smooth closed point of 
a variety over the field $k$. Then in Corollary \ref{constr} a natural, 
continuous in a certain sense homomorphism from the group of cycles of 
codimension $q$ on $X\times_kK$, refining the Abel--Jacobi map:
\begin{equation} \label{Abel-Jacobi}
CH^q(X\times_kK)\longrightarrow\bigoplus_{p=1}^q
gr_F^{q-p}N^{q-p}H^{2q-p}_{dR/k}(X)\otimes_k
\Omega^p_{A^h/k,{\rm closed}},\end{equation}
where $F^{\bullet}$ and $N^{\bullet}$ are the Hodge and the
coniveau filtrations on the de Rham cohomology, respectively, and 
$\Omega^p_{A^h/k}$ is the space of differential $p$-forms on $A^h$. 
The classical Abel--Jacobi map corresponds to the composition 
of this homomorphism with the projection to the summand $p=1$. 
A version of this map has first been studied by V.Srinivas in 
\cite{srinivas}. 

This homomorphism is by no means injective. However, 
its appearance reminds the conjectural formula 
$Ext_{{\cal MM}}^p({\Bbb Q},H^{2q-p}(X)(q))$ 
for the successive quotients of the Chow group $CH^q(X)_{{\Bbb Q}}$ 
with respect to a conjectural filtration on it (cf. \cite{height}). 

Partially motivated by calculations of infinitesimal 
deformations of the Chow groups by Bloch 
(cf. \cite{alg-cycle}, Lecture 6) and Stienstra \cite{stienstra}, 
one could hope its composition with the canonical embedding into the space 
$$\bigoplus_{p=1}^qgr_F^{q-p}N^{q-p}H^{2q-p}_{dR/k}(X)\otimes_k
\lim_{\longleftarrow_M}d(\Omega^{p-1}_{A_M/k}),$$ 
where $A_M=A^h/{\frak m}^M$ and ${\frak m}$ is the maximal ideal, 
is dominant. 

It is shown in Proposition \ref{density} together with the 
concluding remark on p.\pageref{general-mod-lef} that whenever the 
inverse Lefschetz operator 
$H^{2\dim X-q}(X)(\dim X)\stackrel{\sim}{\longrightarrow}H^q(X)(q)$ is 
induced by a self-correspondence on $X$ (e.g., for $q=\dim X$, and 
conjecturally, for arbitrary $q$) the composition is actually dominant. 

I am grateful to Alexander Beilinson for inspiring discussions on 
related topics several years ago and pointing out a mistake in a 
previous version of this note. 

\section{A cycle class $cl_{\bullet}$} 
Let $X$ be a smooth $n$-dimensional projective variety over an 
algebraically closed field $k$ of characteristic 0. Fix an 
auxiliary smooth projective variety $Y$ over 
$k$ with the field of rational functions $K$. 

Using the Poincar\'{e} duality 
$H^{2q}_{dR/k}(X)={\rm Hom}_k(H^{2n-2q}_{dR/k}(X),k)$, one can define 
the class map $cl_{dR}:CH^q(X)\longrightarrow H^{2q}_{dR/k}(X)$ on the 
classes of irreducible subvarieties as 
$[Z]\longmapsto(\omega\mapsto i^{\ast}_Z\omega)$, where 
$i_Z:\widetilde{Z}\longrightarrow X$ is a desingularization 
of $Z$ and $H^{2n-2q}_{dR/k}(\widetilde{Z})$ is canonically identified 
with $k$ via the trace isomorphism. 
\begin{lemma}
Composition of the class map $cl_{dR}$ with the K\"{u}nneth
isomorphism induces a homomorphism $CH^q(X\times_kY)\longrightarrow
\bigoplus_{r=-q}^qH^{q+r}_{dR/k}(X)\otimes_kH^{q-r}_{dR/k}(Y)$. 
that factors through
\begin{equation}
\label{class-ver}
CH^q(X\times_kY)\longrightarrow\bigoplus_{r=0}^qN^rH^{q+r}_{dR/k}(X)\otimes_k
H^{q-r}_{dR/k}(Y)\oplus\bigoplus_{r=1}^qH^{q-r}_{dR/k}(X)\otimes_k
N^rH^{q+r}_{dR/k}(Y).
\end{equation}
\end{lemma}
{\it Proof.} We need to show that for any subvariety 
$Z\subset X\times_kY$ of codimension $q$ and its desingularization 
$\widetilde{Z}$ composition
$$H^{2n-q+r}_{dR/k}(X)\stackrel{{\rm pr}_X^{\ast}}{\longrightarrow}
H^{2n-q+r}_{dR/k}(\widetilde{Z})\stackrel{{\rm pr}_{Y\ast}}
{\longrightarrow}H^{q+r}_{dR/k}(Y)$$
factors through
$H^{2n-q+r}_{dR/k}(X)\longrightarrow N^rH^{q+r}_{dR/k}(Y)$. 
Note, that thanks to the weak Lefschetz theorem, 
$H^{2n-q+r}_{dR/k}(X)=N^{n-q+r}H^{2n-q+r}_{dR/k}(X)$, so our 
composition factors through 
$$H^{2n-q+r}_{dR/k}(X)\stackrel{{\rm pr}_X^{\ast}}{\longrightarrow}
N^{n-q+r}H^{2n-q+r}_{dR/k}(\widetilde{Z}),$$
and, since $Z$ is of relative dimension $n-q$ over $Y$, we have 
$N^{n-q+r}H^{2n-q+r}_{dR/k}(\widetilde{Z})
\stackrel{{\rm pr}_{Y\ast}}{\longrightarrow}N^rH^{q+r}_{dR/k}(Y)$.
\hfill $\Box$

\begin{lemma} \label{class-non-gen-ver}
For any smooth rpojective varieties $X$ and $Y$ over $k$ there is a 
natural homomorphism 
$$CH^q(X\times_kk(Y))\longrightarrow\oplus_{p=0}^q
gr^r_FN^rH^{q+r}_{dR/k}(X)\otimes_k\Gamma(X,\Omega^{q-r}_{Y/k}).$$
\end{lemma}
{\it Proof.} Suppose $D$ is an effective irreducible divisor on $Y$. 
Then for its desingularization $\widetilde{D}$ the composition of the 
inclusion map 
$CH^{q-1}(X\times_k\widetilde{D})\longrightarrow CH^q(X\times_kY)$ with the 
class map (\ref{class-ver}) commutes with composition of the corresponding 
class with $Y$ replaced by $\widetilde{D}$ with the Gysin map and factors 
through the subspace 
$\bigoplus_{r=2-q}^{q-2}H^{q+r}_{dR/k}(X)\otimes_kN^1H^{q-r}_{dR/k}(Y)$, 
i.e., the following diagram commutes 
$$
\begin{diagram}
\node{CH^{q-1}(X\times_k\widetilde{D})} \arrow[2]{e} \arrow{s} \node{} 
\node{CH^q(X\times_kY)} \arrow{s} \\
\node{H^{2(q-1)}_{dR/k}(X\times_k\widetilde{D})} 
\arrow{e,t}{id\otimes{\rm Gysin}}  
\node{\bigoplus_{r=2-q}^{q-2}H^{q+r}_{dR/k}(X)\otimes_kN^1H^{q-r}_{dR/k}(Y)} 
\arrow{e,t}{\subseteq} \node{H^{2q}_{dR/k}(X\times_kY)} 
\end{diagram}
$$ 

Since one has the exact localization sequence 
$$\bigoplus_{D:{\rm divisors~on}~Y}CH^q(X\times_kD)\longrightarrow 
CH^q(X\times_kY)\longrightarrow CH^q(X\times_kk(Y))\longrightarrow 0,$$ 
the homomorphism (\ref{class-ver}) leads to a well-defined 
homomorphism from the Chow group \\
$CH^q(X\times_kk(Y))$ to the quotient of the 
right hand side of (\ref{class-ver}) by the subspace \\
$\bigoplus_{r=1}^qH^{q-r}_{dR/k}(X)\otimes_kN^1H^{q+r}_{dR/k}(Y)$, 
i.e., 
$$CH^q(X\times_kk(Y))\longrightarrow\bigoplus_{r=0}^qN^r
H^{q+r}_{dR/k}(X)\otimes_k\frac{H^{q-r}_{dR/k}(Y)}{N^1H^{q-r}_{dR/k}(Y)}.$$ 
In fact, the space $H^{q-r}_{dR/k}(Y)/N^1H^{q-r}_{dR/k}(Y)$ is 
a birational invariant of $Y$ and can be considered as a subspace of 
$H^{q-r}_{dR}(k(Y)/k)$, or as a subspace of $\Omega^{q-r}_{k(Y)/k}$. 

Also, from the K\"{u}nneth decomposition 
$$F^qH^p_{dR/k}(X\times_kY)=
\sum_{s=0}^q\sum_{t=0}^pF^sH^t_{dR/k}(X)\otimes_kF^{q-s}H^{p-t}_{dR/k}(Y)$$ 
one sees that any class in $F^q$ and simultaneously in the
space 
$$F^rH^{q+r}_{dR/k}(X)\otimes_kH^{q-r}_{dR/k}(Y)
\quad\mbox{\rm but not in the subspace}\quad
F^{r+1}H^{q+r}_{dR/k}(X)\otimes_kH^{q-r}_{dR/k}(Y)$$
belongs to $F^rH^{q+r}_{dR/k}(X)\otimes_kF^{q-r}H^{q-r}_{dR/k}(Y)$.
Therefore, we have a homomorphism 
$$CH^q(X\times_kk(Y))\longrightarrow\bigoplus_{r=0}^qgr_F^rN^r
H^{q+r}_{dR/k}(X)\otimes_k\Gamma(Y,\Omega^{q-r}_{Y/k}).\quad\quad\Box$$

At this point, we may replace the field of rational functions on $Y$
with an arbitrary field extension $K\subset k$ to obtain the 
following class map 
\begin{equation}
\label{class}
cl_q:CH^q(X\times_kK)\longrightarrow
\bigoplus_{r=0}^qgr_F^rN^rH^{q+r}_{dR/k}(X)\otimes_k\Omega^{q-r}_{K/k}.
\end{equation}

\section{Invariants of cycles over fraction fields of henselizations}
For any pair of rings $B\subseteq C$ denote by 
$\Omega^{\bullet}_{C/B}$ the exterior $C$-algebra of the module 
$\Omega^1_{C/B}$ of differentials on $C$ relative to $B$. 
\begin{lemma} \label{integrality}
Let $k$ be a field of characteristic zero and 
$A$ be a local ring of a closed regular point of a variety over
$k$ with the fraction field $K$. Let $X$ be a smooth projective 
variety and $\sigma:k(X)\hookrightarrow K$ an embedding of fields. 
Obviously, one can extend $\sigma:k(X)\hookrightarrow K$ to a homomorphism 
$\Omega^{\bullet}_{k(X)/k}\hookrightarrow\Omega^{\bullet}_{K/k}$ of 
differential $k$-algebras in a unique way.
Then the restriction of this homomorphism 
$\Gamma(X,\Omega^{\bullet}_{X/k})\longrightarrow\Omega^{\bullet}_{K/k}$ 
factors through an embedding 
$\Gamma(X,\Omega^{\bullet}_{X/k})\longrightarrow\Omega^{\bullet}_{A/k}
\subset\Omega^{\bullet}_{K/k}$. 
\end{lemma}
{\it Proof.} Since $K=k(Y)$ is the function field of a smooth 
projective variety $Y$, we may suppose that $\sigma$ is induced 
by a morphism  $Y\longrightarrow\!\!\!\!\!\rightarrow X$. Then we may 
find a smooth projective variety $Y'$ birational to $Y$ with $A$ as 
the local ring of $Y'$ at a (closed) point $y$. We are done, since 
$$\sigma:\Gamma(X,\Omega^{\bullet}_{X/k})\longrightarrow
\Gamma(Y,\Omega^{\bullet}_{Y/k})=\Gamma(Y',\Omega^{\bullet}_{Y'/k})
\subset\Omega^{\bullet}_{{\cal O}_{Y',y}/k}.
\qquad\qquad\qquad\qquad\qquad\Box$$ 

This lemma together with the homomorphism $cl_q$ from 
(\ref{class}) enables us to define the homomorphism 
(\ref{Abel-Jacobi}). 
\begin{cor} \label{constr} For any smooth projective variety $X$ over 
$k$ and the field of fractions $K$ of $A^h$ the map 
$$\bigoplus_{q\geq 0}CH^q(X\times_kK)\longrightarrow
\bigoplus_{0\leq p\leq q}gr^{q-p}_FN^{q-p}H^{2q-p}_{dR/k}(X)
\otimes_k\Omega^p_{A^h/k,{\rm closed}}$$
is a homomorphism of the graded rings. \hfill $\Box$
\end{cor}
\begin{lemma} \label{reduction-1}
Suppose that for a local domain $A\supset k$ with the 
fraction field $K$, an integer $q\geq 0$, 
for any integer $n$ and any smooth $n$-dimensional 
projective variety $X$ over $k$ and any integer $M>0$ 
$$CH_q(X\times_kK)\longrightarrow\!\!\!\!\!\rightarrow 
H^{n-q}(X,{\cal O}_X)\otimes_kd(\Omega^{n-q-1}_{A_M/k}),$$
where $A_M=A/{\frak m}_A^M$. Then the map 
$$CH_q(X\times_kK)\longrightarrow\!\!\!\!\!\rightarrow\bigoplus_{p=1}^{n-q}
gr_F^{n-p-q}N^{n-p-q}H^{2(n-q)-p}_{dR/k}(X)\otimes_kd(\Omega^{p-1}_{A_M/k})$$
is surjective for any integer $M>0$. 
\end{lemma}
{\it Proof.} We proceed by induction on dimension of $X$. The case 
$\dim X\leq q+1$ is trivial. 

Denote by $\widetilde{D}$ a desingularization of divisor $D$ on $X$. 
Then Lemma follows from the commutativity of the diagram 
$$
\begin{diagram}
\node{CH_q(X\times_kK)} \arrow[2]{e} \node{} \node{\bigoplus_{p=1}^{n-q}
gr_F^{n-p-q}N^{n-p-q}H^{2(n-q)-p}_{dR/k}(X)\otimes_kd(\Omega^{p-1}_{A_M/k})} \\
\node{\bigoplus_{D\in X^1}CH_q(\widetilde{D}\times_kK)} 
\arrow[2]{e,t}{{\rm onto}} \arrow{n} \node{} 
\node{\bigoplus_{p=1}^{n-q-1}gr_F^{n-p-q-1}N^{n-p-q-1}H^{2(n-q-1)-p}_{dR/k}
(\widetilde{D})\otimes_kd(\Omega^{p-1}_{A_M/k})} 
\arrow{n,r}{{\rm Gysin}\otimes id}
\end{diagram}
$$
where the surjectivity of the Gysin maps 
$\bigoplus_{D\in X^1}N^{\ast-1}H^{\bullet-2}(\widetilde{D})
\longrightarrow\!\!\!\!\!\rightarrow N^{\ast}H^{\bullet}(X)$ 
should be used.\footnote{To check the surjectivity of 
the Gysin maps one uses alternative description of the coniveau filtration 
due to Grothendieck: $N^aH^b(X)={\rm Im}(\oplus_{Z\in X^a}
H^{b-2a}(\widetilde{Z})\stackrel{{\rm Gysin}}{\longrightarrow}H^b(X))$, 
cf. \cite{brau}, formula (10.7) or (9.17), and a correction in footnote 
on p.300 of \cite{hodge}.} \hfill $\Box$ 

\begin{lemma} \label{curves}
Let $R$ be an integral domain, and $\|f_{jk}\|$ 
be an invertible $g\times g$-matrix with entries in the algebra 
$R[[t]]$ of power series in one variable $t$. 

Then there exists a unique collection of power series 
$\phi_l\in tR[[t]]$ for $1\leq l\leq g$ such that 
$$\sum_{j=1}^gf_{ij}(\phi_j)\frac{d\phi_j(t)}{dt}=\left\{
\begin{array}{ll} 1 & \mbox{if $i=1$,} \\
0 & \mbox{otherwise.} \end{array} \right.$$
\end{lemma}
{\it Proof.} Let $\phi_l=\phi_l(t)\in tR[[t]]$. 
We need to solve the system 
$$\left\{\begin{array}{c}\sum_{l=1}^gf_{1l}(\phi_l(t))\phi'_l(t)=1 \\
\sum_{l=1}^gf_{2l}(\phi_l(t))\phi_l'(t)=0 \\
\dots \\
\sum_{l=1}^gf_{gl}(\phi_l(t))\phi_l'(t)=0 \end{array} \right.$$
and claim that there is a unique solution of this system. 
This is equivalent to saying that there is a unique collection 
of $\phi^{(j)}_l(0)$. We prove the latter by induction on $j$, the 
case $j=0$ being trivial. 
Taking $j$th derivative at $t=0$ leads to the system 
\begin{equation} \label{system-j}
\left\{\begin{array}{c} 
\sum_{l=1}^gf_{1l}(0)\phi^{(j+1)}_l(0)=\langle\mbox{polynomial
in derivatives of $\phi$'s at 0 of orders $\leq j$}\rangle \\
\sum_{l=1}^gf_{2l}(0)\phi_l^{(j+1)}(0)=\langle\mbox{polynomial
in derivatives of $\phi$'s at 0 of orders $\leq j$}\rangle \\
\dots \\
\sum_{l=1}^gf_{gl}(0)\phi_l^{(j+1)}(0)=\langle\mbox{polynomial
in derivatives of $\phi$'s at 0 of orders $\leq j$}\rangle
\end{array} \right. \end{equation}

Since the $g\times g$-matrix $\|f_{jl}\|$ is invertible, 
the system (\ref{system-j}) has a unique solution \\
$(\phi_1^{(j+1)}(0),\dots,\phi_g^{(j+1)}(0))$. \hfill $\Box$

\begin{lemma} \label{points-choice}
Let $U$ be a variety over an algebraically closed field $k$ and 
$f_1,\dots,f_g\in{\cal O}(U)$ be a collection of regular
functions on $U$, linear independent over the field $k$. 

Then there exist such a collection of points $p_1,\dots,p_g$ 
that the $g\times g$-matrix $\|f_j(p_k)\|$ is invertible together
with its first-row minors. 
\end{lemma}
{\it Proof.} We proceed by induction on $g$, the case $g=1$ being
trivial. 

The matrix $\|f_j(p_k)\|$ is not invertible for any
collection of points $p_1,\dots,p_g$ means that 
the function 
$$\det\|f_j(p_k)\|:\underbrace{U\times\dots\times U}_
{g~{\rm times}}\longrightarrow k$$
is identically zero. 

Decomposing the determinant via the first row, we get 
$$\det\|f_k(p_l)\|=\sum_{j=1}^g(-1)^{j-1}
\det\|f_k(p_l)\|_{k\neq j,l>1}f_j(p_1),$$ 
where, under the induction assumption, 
$\det\|f_k(p_l)\|_{k\neq j,l>1}$ is 
a non-zero function in $p_2,\dots,p_g$. 

After fixing sufficiently general $p_2,\dots,p_g$ 
we get $\det\|f_j(p_k)\|=\sum_{j=1}^ga_jf_j(p_1)$, where $a_j$ 
are non-zero elements of $k$, and therefore, we get the $k$-linear
dependence of the functions $f_1,\dots,f_g$. \hfill $\Box$

\begin{prop} \label{density}
Let $k$ be a field of characteristic zero, $N>0$ an integer, 
$K$ the algebraic closure of the field of rational functions 
in the fraction field of $k[[t_1,\dots,t_N]]$. 

Then for any smooth $n$-dimensional projective variety $X$ over $k$ 
and an integer $M>0$ there is a natural surjection 
\begin{equation} \label{infinites} 
CH_0(X\times_kK)\longrightarrow\!\!\!\!\!\rightarrow\bigoplus_{q=1}^n
H^n(X,\Omega^{n-q}_{X/k})\otimes_kd(\Omega^{q-1}_{A_M/k}),
\end{equation}
where $A_M=k[t_1,\dots,t_N]/(t_1,\dots,t_N)^M$.
\end{prop}
{\it Proof.} Since surjectivity is preserved by taking Galois
invariants, we may suppose that the field $k$ is algebraically closed.
By Lemma \ref{reduction-1} it is enough to check that 
the composition of (\ref{infinites}) with the projection to 
$H^n(X,{\cal O}_X)\otimes_kd(\Omega^{n-1}_{A_M/k})$ is surjective. 

We fix a basis $\{\omega_0,\dots,\omega_N\}$ of the space
$\Gamma(X,\Omega^n_{X/k})$ and choose functions $x_1,\dots,x_n$
algebraically independent over $k$. Set $f_j=\omega_j/\omega_0$. 

Fix an open subset $U$ where $(x_1,\dots,x_n):U\longrightarrow{\Bbb A}^n_k$ 
is an \'{e}tale morphism and the functions $f_j$ are regular. 
Choose points $p_1,\dots,p_g$ as in Lemma \ref{points-choice}. 
This guarantees that the matrix $\|f_i(p_j)\|$ is invertible. 

Then for each point $p_j$ we embed the completion of 
${\cal O}(U)$ at $p_j$ into the algebra $k[[t_1,\dots,t_n]]$ 
by sending $x_l$ to $x_l(p_j)+t_l$. 
Denote by $f_{ij}$ the image of $f_i$ under this embedding. 

Note, that the matrix $\|f_{jl}\|$ with entries in 
$k[[t_1,\dots,t_n]]$ is invertible together with its first-row
minors, since its reduction modulo the maximal ideal in 
$k[[t_1,\dots,t_n]]$ is the matrix $\|f_i(p_j)\|$. 

Finally, we consider $f_{jl}$ as elements of $R[[t_1]]$ with 
$R=k[[t_2,\dots,t_n]]$, and apply Lemma \ref{curves} to show 
that there exists a collection of non-zero 
($\partial\phi_l/\partial t_1(0,\dots,0)\neq 0$) formal series 
$\phi_l(t_1,\dots,t_n)\in t_1k[[t_1,\dots,t_n]]$ such that 
$$\sum_lf_{il}(\phi_l(t_1,\dots,t_n),t_2,\dots,t_n)
d\phi_l\wedge dt_2\wedge\dots\wedge dt_n
=\left\{\begin{array}{ll} dt_1\wedge\dots\wedge dt_n & \mbox{if $i=1$,} \\
0 & \mbox{otherwise.}\Box \end{array} \right.$$

{\sc Remark}. \label{general-mod-lef} It follows from the standard 
conjecture on algebraicity of the inverse Lefschetz operator 
$H^{2\dim X-q}(X)(\dim X)\stackrel{\sim}{\longrightarrow}H^q(X)(q)$ 
that the map (\ref{Abel-Jacobi}) is dominant. 

{\it Proof.} Thanks to Lemma \ref{reduction-1} we only need to
check that for any element $\beta$ of the $k$-vector space 
$H^q(X,{\cal O}_X)\otimes_kd(\Omega^{q-1}_{A_M/k})$
there exists a cycle $\alpha$ in $CH^q(X\times_kK)$ with the
image $\beta$ under the composition of the map
(\ref{Abel-Jacobi}) with the corresponding projection. 

Fix a smooth $q$-dimensional plane section $W$ of $X$. 
Denote by $\overline{\beta}$ the image of $\beta$ under the
isomorphism 
$$H^q(X,{\cal O}_X)\otimes_kd(\Omega^{q-1}_{A_M/k})
\stackrel{\cup[W]\otimes id}{\longrightarrow} 
H^n(X,\Omega^{n-q}_{X/k})\otimes_kd(\Omega^{q-1}_{A_M/k}).$$

Due to Proposition \ref{density}, there is a 0-cycle
$\overline{\alpha}$ with the image $\overline{\beta}$ under the map 
(\ref{Abel-Jacobi})$=$(\ref{infinites}). Suppose that
$\overline{\alpha}$ is, in fact, defined over a field
$k(Y)\subset K$ finitely generated over $k$, and $Y$ is its
smooth projective model over $k$. 
Let $\widetilde{\alpha}\in CH^n(X\times_kY)$ correspond to the 0-cycle 
$\overline{\alpha}\in CH_0(X\times_kK)$. Assuming that the inverse Lefschetz 
operator $H^{2n-q}(X)(n)\stackrel{\sim}{\longrightarrow}H^q(X)(q)$ is
induced by a correspondence $[L]\in CH^q(X\times_kX)$, we set 
$\widehat{\alpha}=[L\times\Delta_Y](\widetilde{\alpha})$.

The image of the cycle $\widehat{\alpha}$ in $CH^q(X\times_kK)$ 
gives us the desired cycle $\alpha$. \hfill $\Box$

\end{document}